\documentclass{amsart}
\usepackage{mathbbol}
\usepackage{amsmath, amsthm, amssymb}
\usepackage{txfonts}

\usepackage[usenames]{color}

\newtheorem*{thm}{Theorem}
\newcommand{\bt}{\begin{thm}}
\newcommand{\et}{\end{thm}}

\newtheorem{cor}{Corollary}
\newcommand{\bc}{\begin{cor}}
\newcommand{\ec}{\end{cor}}

\newtheorem*{lem}{Lemma}
\newcommand{\bl}{\begin{lem}}
\newcommand{\el}{\end{lem}}

\newtheorem*{prop}{Proposition}
\newcommand{\bp}{\begin{prop}}
\newcommand{\ep}{\end{prop}}

\newtheorem{defn}{Definition}
\newcommand{\bd}{\begin{defn}}      
\newcommand{\ed}{\end{defn}}

\newtheorem{rmrk}{Remark}
\newcommand{\br}{\begin{rmrk}}
\newcommand{\er}{\end{rmrk}}

\newcommand{\secref}[1]{Section~\ref{#1}}

\newcommand{\R}{\mathbb{R}}
\newcommand{\Z}{\mathbb{Z}}

\newcommand{\diam}{\operatorname{diam}}

\newcommand{\systole}{\operatorname{sys}_1}

\newcommand{\spt}{\operatorname{spt}}

\newcommand{\Vol}{\operatorname{Vol}}
\newcommand{\fillvol}{{\operatorname{Fillvol}}}

\newcommand{\bdry}{\partial}

\begin{document}

\title[Gromov's filling inequality]{A short proof of Gromov's filling inequality}

\author{Stefan Wenger}

\address
  {Courant Institute of Mathematical Sciences\\
   251 Mercer Street\\
   New York, NY 10012}
\email{wenger@cims.nyu.edu}

\date{March 29, 2007}



\begin{abstract}
 We give a very short and rather elementary proof of Gromov's filling volume inequality for $n$-dimensional Lipschitz cycles (with integer and $\Z_2$-coefficients) in $L^\infty$-spaces.
 This inequality is used in the proof of Gromov's systolic inequality for closed aspherical Riemannian manifolds and is often regarded as the difficult step therein.
\end{abstract}

\maketitle


\bigskip

\section{Introduction}

Gromov's well-known systolic inequality in \cite{Gromov-filling} asserts that for any $n$-dimensional closed aspherical Riemannian manifold $(M^n,g)$,
the shortest length $\systole(M^n,g)$ of a non-contract\-ible loop in $M^n$ satisfies
\begin{equation}\label{equation:systolic}
 \systole(M^n,g)\leq c_n \Vol(M^n,g)^\frac{1}{n},
\end{equation}
with $c_n$ only depending on $n$. The proof of this relies on the following isoperimetric inequality, often regarded as
the difficult step in the proof.

\bt[{Gromov \cite[2.3]{Gromov-filling}}]
 Let $X$ be an $L^\infty$-space and $n\geq 1$.
 Then the filling volume of any $n$-dimensional singular Lipschitz cycle $z$ in $X$ with integer or $\Z_2$-coefficients satisfies
 \begin{equation*}
  \fillvol(z)\leq C_n \Vol(z)^{1+\frac{1}{n}},
 \end{equation*}
 where $C_n$ depends only on $n$.
\et

Here, $X$ is said to be an $L^\infty$-space if $X=L^\infty(\Omega)$, the space of bounded functions on some set $\Omega$, endowed with the supremum norm.
Furthermore, $\fillvol(z)$ is the least volume of an $(n+1)$-dimensional singular Lipschitz chain in $X$ with boundary $z$. 

The proof given by Gromov in \cite{Gromov-filling} in fact applies not only to $L^\infty$ but to arbitrary Banach spaces and to Riemannian 
manifolds admitting cone type inequalities and also works for cycles with real coefficients. For the systolic inequality, however, the above theorem is enough.
In \cite{Wenger-GAFA} it was later shown, that Gromov's result (for cycles with integer coefficients) extends to the setting of complete metric spaces admitting cone type inequalities. 
The arguments in \cite{Wenger-GAFA} rely on the theory of integral currents in metric spaces \cite{Ambrosio-Kirchheim}; these give a generalization of Lipschitz chains 
(only) with integer coefficients.

The purpose of this note is to provide a short and rather elementary proof of the above theorem. The main idea stems from \cite{Wenger-GAFA}, but
our proof does not use at all the theory of currents. 
Apart from being less than two pages, it enjoys the following features different from \cite{Gromov-filling}: 
Firstly, it works also (and without modification) for ${\rm CAT}(0)$-spaces, i.e.~simply connected geodesic metric spaces of non-positive curvature in the sense of Alexandrov. Secondly, 
it does not use the Federer-Fleming isoperimetric inequality in $\R^N$ and in fact yields a new proof for the isoperimetric inequality in Euclidean space.

\bigskip

{\bf Acknowledgments:} The author thanks Misha Katz for encouraging him to write up this note. The paper was written while the author was visiting the Max Planck Institute
of Mathematics in Bonn.

\section{Singular Lipschitz chains}\label{section:notion-volume}
Recall from \cite{Gromov-filling} that a singular Lipschitz $n$-chain is a formal finite sum $c = \sum_{i=1}^k m_i\varphi_i,$
where the $\varphi_i$ are Lipschitz maps $\varphi_i: \Delta^n\to X$ and the $m_i$ values in $\Z$ or $\Z_2$. Here $\Delta^n\subset\R^{n+1}$ is the standard $n$-simplex. 
The boundary $\bdry c$ is by definition the singular Lipschitz chain $\bdry c = \sum_{i=1}^k m_i\bdry\varphi_i$, where $\bdry\varphi_i$ is the $(n-1)$-chain induced by 
restricting $\varphi_i$ to the faces of $\Delta^n$.
The volume of $c$ is $\Vol(c):= \sum|m_i|\Vol(\varphi_i)$, where $\Vol(\varphi_i)$ is the `parametrized' volume of $\varphi$. In the setting of Banach spaces there exists more than one
suitable definition of volume. Therefore, $\Vol(\varphi_i)$ may be chosen to be the parametrized Hausdorff, the Holmes-Thompson, or Gromov's mass$*$ volume.
Since all these volumes agree up to a factor at most $n^{n/2}$ the choice of volume is irrelevant for the proof of the isoperimetric inequality.
In what follows singular Lipschitz $n$-chains will simply be called $n$-chains.

For the proof of the isoperimetric inequality the following three basic facts will be needed: 
There exist constants $A_n\geq 2^{-1}n^{-n/2}\omega_n$ and $D_n, E_n$ depending only on $n$ and the chosen volume such that the following properties hold
($\omega_n$ denotes the volume of the unit ball in $\R^n$):
Let $X$ be a metric space and $z$ an $n$-cycle in $X$.\\

{\bf Coarea formula:} Given a Lipschitz function $\varrho: X\to\R$ with Lipschitz constant $1$ then
 \begin{equation}\label{equation:coarea-inequality}
  \Vol_{n-1}(z\cap \{\varrho=r\})\leq E_n\frac{d}{dr}\Vol_n(z\cap\{\varrho\leq r\})
 \end{equation}
 for almost every $r\in\R$. Furthermore, $z\cap\{\varrho=r\}$ is an $(n-1)$-cycle for almost every $r\in\R$ (after a minor `smoothening' of $\{\varrho=r\}$).\\
 
{\bf Cone inequality:} If $X$ is a ${\rm CAT}(0)$ or a Banach space and $x_0\in X$ then
the $(n+1)$-chain $c$ in $X$ obtained by joining each point in $z$ by the geodesic line (respectively, straight line
if $X$ is a Banach space) with $x_0$ satisfies $\bdry c=z$ and
\begin{equation}\label{equation:cone-inequality}
 \Vol_{n+1}(c)\leq D_nR\Vol_n(z),
\end{equation}
where $R$ is the smallest number such that $z$ is contained in the ball or radius $R$ around $x_0$.\\

{\bf Lower density estimate:} For almost every $x$ in the support of $z$
 \begin{equation}\label{equation:lower-density}
  \liminf_{r\to 0^+} \frac{1}{r^n}\Vol_n(z\cap B(x,r))>A_n.
 \end{equation}

If $X$ is a Banach space and $\Vol$ is Gromov's mass$*$ volume then $D_n\leq 1$ and $E_n=1$. 
If $X$ is a ${\rm CAT}(0)$-space then all the definitions of volume agree and we have $D_n=\frac{1}{n+1}$ and $E_n=1$.

Gromov's proof \cite{Gromov-filling} uses the coarea formula and the cone inequality, but not the lower density estimate 
(which is false for chains with real coefficients).
Note that \eqref{equation:lower-density} is clear when $z$ is a piecewise $C^1$-submanifold. In our arguments this will be enough (by a `smoothening' argument). 
For Lipschitz cycles the proof is given in \cite{Kirchheim}.

\section{The short proof}\label{section:Gromov-theorem-proof}
Fix a definition of volume, let $n\geq 1$ and let $A_n, D_n, E_n$ be the corresponding values from \secref{section:notion-volume}. We abbreviate 
$\Vol=\Vol_n$.
Let $0<\varepsilon<1$ be given.
%

\bp
 Let $X$ be a metric space and let $z$ be an $n$-cycle in $X$, $n\geq 2$. There exist finitely many pairwise disjoint closed balls $B_i\subset X$, $i=1,\dots, k$,
 with the following properties:
 \begin{enumerate}
  \item The volume of $z$ contained in each ball is `not too small':  
   \begin{equation*}
    \Vol(z\cap B_i)\geq 4^{-n}A_n\varepsilon\diam(B_i)^n.
   \end{equation*}
  \item The restriction $z\cap B_i$ is an $n$-chain whose boundary has `small' volume:
   \begin{equation*}
    \Vol_{n-1}(\bdry (z\cap B_i))\leq E_n(A_n\varepsilon)^{\frac{1}{n}}n \Vol(z\cap B_i)^{\frac{n-1}{n}}.
   \end{equation*}
  \item An essential part of the volume of $z$ is contained in the union of these balls: 
   \begin{equation*}
    \Vol\left(z\cap \bigcup B_i\right)\geq \frac{1}{5^n}\Vol(z).
   \end{equation*}
 \end{enumerate}
\ep

\begin{proof}
 For $x\in\spt z$ and $r\geq 0$ define $V(x,r):= \Vol(z\cap B(x,r))$ and
 \begin{equation*}
  r_0(x):= \max\{r\geq 0: V(x,r)\geq A_n\varepsilon r^n\}.
 \end{equation*}
 Note that $0<r_0(x)<\infty$ for almost every $x\in\spt z$ by the lower density estimate; moreover 
 \begin{equation*}
  V(x,5r_0(x))< 5^nA_n\varepsilon r_0(x)^n= 5^n V(x,r_0(x)).
 \end{equation*}
 By Vitali's $5r$-covering lemma there exist finitely many points $x_i\in \spt z$, $i=1,\dots, k$, such that the balls $B(x_i, 2r_0(x_i))$ are pairwise disjoint, the 
 balls $B(x_i, 5r_0(x_i))$ cover $\spt z$ and
 \begin{equation*}
  \Vol\left(z\cap \bigcup B(x_i, r_0(x_i))\right)\geq \frac{1}{5^n}\Vol(z). 
 \end{equation*}
 Fix $i\in\{1,\dots,k\}$.  By the the definition of $r_0(x_i)$ there exists a non-negligible set of points $r\in(r_0(x_i), 2r_0(x_i))$ with 
  \begin{equation*}
   \frac{d}{dr}V(x_i,r)< (A_n\varepsilon)^\frac{1}{n}n V(x_i,r)^{\frac{n-1}{n}},
  \end{equation*}
  therefore, by the coarea inequality,
  \begin{equation}\label{equation:bdryvol}
    \Vol_{n-1}(\bdry(z\cap B(x_i,r))=\Vol_{n-1}(z\cap\{x: d(x,x_i)=r\})< E_n(A_n\varepsilon)^{\frac{1}{n}}n  \Vol(z\cap B(x_i, r))^{\frac{n-1}{n}}.
  \end{equation}
  Choose an $r$ such that \eqref{equation:bdryvol} holds and set $B_i:=B(x_i, r)$. The so defined $B_i$ clearly satisfy (i), (ii) and (iii).
 \end{proof}

Let now $X$ be a ${\rm CAT}(0)$ or an $L^\infty$-space.

\begin{proof}[{Proof of Gromov's theorem}]
 The proof is by induction on $n$ and the case $n=1$ is trivial, since the diameter of a closed curve is bounded by its length and thus the isoperimetric inequality 
 is a direct consequence of the cone inequality. 
 Suppose now that $n\geq 2$ and that the statement of the theorem holds for $(n-1)$-cycles with some constant $C_{n-1}\geq 1$.
 Set 
 \begin{equation*}
  \varepsilon:=\min\left\{\frac{1}{4^{n-1}C_{n-1}^{n-1}A_nE^n_nn^n},  \frac{1}{2}\right\},
 \end{equation*}
 let $z$ be an $n$-cycle in $X$ and choose a ball $B$ of finite radius which contains $z$.
 Let $B_1,\dots, B_k$ be balls as in the proposition. By the isoperimetric inequality in dimension $n-1$ we can choose for each $i=1,\dots, k$ an $n$-chain $c_i$ 
 satisfying $\bdry c_i=\bdry (z\cap B_i)$ and
 \begin{equation}\label{equation:chain-upperbound-mass}
  \Vol(c_i)\leq C_{n-1} \Vol_{n-1}(\bdry(z\cap B_i))^{\frac{n}{n-1}} 
  \leq\frac{1}{4}\Vol(z\cap B_i).
 \end{equation}
 Here the second inequality follows from (ii) and the definition of $\varepsilon$.
 We may of course assume that $c_i$ is contained in $B_i$ since otherwise we can project it to $B_i$ via a $1$-Lipschitz projection $P:X\to B_i$ (and this decreases the volume). 
 If $X=L^\infty(\Omega)$ then $P(f)(y):= \operatorname{sgn}(f(y))\min\{|f(y)|, 1\}$. If  $X$ is a ${\rm CAT}(0)$-space then $P$ is the orthogonal projection onto $B_i$.
 Set $\hat{z}_i:= (z\cap B_i) - c_i$ and $$z':= z - \sum_{i=1}^k \hat{z}_i = (z\cap (\cup B_i)^c) - \sum_{i=1}^k c_i.$$
 Observe that these are $n$-cycles and that, 
 by \eqref{equation:chain-upperbound-mass},
 \begin{equation}\label{equation:mass-ball}
  \frac{3}{4}\Vol(z\cap B_i)\leq \Vol(\hat{z}_i)\leq \frac{5}{4}\Vol(z\cap B_i).
 \end{equation}
 From the proposition and from \eqref{equation:chain-upperbound-mass}, \eqref{equation:mass-ball} we conclude
 \begin{equation}\label{equation:round}
  \diam(\hat{z}_i)\leq \diam B_i\leq \frac{4}{(A_n\varepsilon)^{\frac{1}{n}}}\Vol(z\cap B_i)^\frac{1}{n} \leq \left(\frac{4^{n+1}}{3A_n\varepsilon}\right)^{\frac{1}{n}}\Vol(\hat{z}_i)^{\frac{1}{n}}
 \end{equation}
 and
 \begin{equation}\label{equation:volume-sum}
  \frac{3}{5}\left[\sum_{i=1}^k\Vol(\hat{z}_i)\right] + \Vol(z')\leq \Vol(z)
 \end{equation}
 as well as $\Vol(z')\leq (1-\frac{3}{4}5^{-n})\Vol(z).$
After further decomposing $z'$ as above we may assume that $\Vol(z')$ is as small as we want.
Let $\hat{c}_i$ and $c'$ be $(n+1)$-chains with boundaries $\hat{z}_i$ and $z'$, respectively, and which satisfy the cone inequality. Note that we can make $\Vol(c')$ 
as small as we wish because $\diam z'\leq \diam B<\infty$.
The $(n+1)$-chain $c:= \hat{c}_1+\dots+\hat{c}_k + c'$ has boundary $z$ and satisfies
 \begin{equation*}
    \Vol(c)\leq \Vol(\hat{c}_1)+\dots+\Vol(\hat{c}_k)+\Vol(c')\leq C_n\Vol(z)^{1+\frac{1}{n}}
 \end{equation*}
 for some $C_n$ only depending on $n$.
 This is a consequence of \eqref{equation:round}, \eqref{equation:volume-sum}, and the fact that 
 for $a_1,\dots,a_k\geq 0$ and $\alpha\geq 1$
  \begin{equation*}
   a_1^\alpha+\dots+a_k^\alpha\leq (a_1+\dots+a_k)^\alpha.
  \end{equation*}
 This concludes the proof.
\end{proof}

By evaluating the expressions in the proof one sees that the constant $C_n$ is bounded above by
\begin{equation*}
 C_n \leq 27nD_nE_nC_{n-1}^{\frac{n-1}{n}}.
\end{equation*}
In particular, if $\Vol$ is Gromov's mass$*$ volume, the isoperimetric constant is bounded from above by $27^nn!$. Note that this is 
smaller than the constant which appears in \cite{Gromov-filling}.

Clearly, our proof works for all metric spaces with cone inequality for which there exist $1$-Lipschitz projections onto closed balls.

\end{document}